\documentclass{amsart}
 \newtheorem{theorem}{Theorem}[section]
       \newtheorem{proposition}[theorem]{Proposition}
       \newtheorem{corollary}[theorem]{Corollary}
       \newtheorem{lemma}[theorem]{Lemma}
       
\theoremstyle{definition}

\newcommand{\la}{\lambda}
\newcommand{\eps}{\varepsilon}
\newcommand{\RR}{\mathbb{R}}
\newcommand{\CC}{\mathbb{C}}
\newcommand{\mU}{{\bf U}}
\newcommand{\mV}{{\bf V}}
\newcommand{\mA}{{\bf A}}

\numberwithin{equation}{section}
\date{July 8, 2006 (file {\tt \jobname.tex} printed \today)}
\title[Empirical measures of weighted sums]{Asymptotic results for empirical measures of weighted sums
of independent random variables}
\author{Bernard Bercu}
\address{
Laboratoire de Statistique et Probabilit\'es, UMR C5583,
Universit\'e Paul Sabatier, 118 Route de Narbonne, 31062 Toulouse
Cedex, France} \email{bercu@cict.fr}

\author{W{\l}odzimierz  Bryc}
\thanks{\noindent Research partially supported by NSF
grant \#DMS-0504198.}
\address{
Department of Mathematical Sciences, University of Cincinnati, 2855
Campus Way, PO Box 210025, Cincinnati, OH 45221-0025, USA}
\email{Wlodzimierz.Bryc@UC.edu}


\begin{document}
\maketitle
\begin{abstract}
We prove that  if a rectangular $r\times n$ matrix with uniformly
small entries and approximately orthogonal rows is applied to the
independent standardized random variables with uniformly bounded
third moments, then the empirical CDF of the resulting partial sums
converges to the  normal CDF with probability one. This implies
almost sure convergence of empirical periodograms, almost sure
convergence of spectra  of circulant and reverse circulant matrices,
and almost sure convergence of the CDF's generated from independent
random variables by independent random orthogonal matrices.

For special trigonometric matrices, the speed of the almost sure
convergence is described by the normal approximation and by the
large deviation principle.
\end{abstract}

\section{Results}
The study of spectra of circulant matrices  lead Massey, Miller and
Shinsheimer \cite{Miller05}  to an almost sure Central Limit Theorem
(CLT) which takes a different form than the celebrated almost sure
CLT discovered by Brosamler \cite{Brosamler88} and Schatte
\cite{Schatte88}. Our main result extends \cite[Theorem
5.1]{Miller05} in several ways: we allow more general weights, we do
not assume identical distributions, we do not assume higher moments
than three,
 and we prove  multivariate convergence. We
consider weighted sums of independent random variables with the
weights that come from matrices with "almost orthogonal" rows and
uniformly small entries. Namely, let $\mU^{(n)}=[u_{i,j}^{(n)}]$ be
a family of real
 $r_n\times n$ matrices, $n\geq 1$. We assume that for some
constants $C, \delta>0$ which do not depend on $n=1,2,\dots$, we
have
\begin{eqnarray}
 \max_{1\leq j\leq n,1\leq k\leq r_n} |u_{k,j}^{(n)}| &\leq& \frac{C}{(\log (1+r_n))^{1+\delta}}, \label{i}
\\
\max_{1\leq k_1,k_2\leq r_n} \left|\sum_{j=1}^n u_{k_1,j}^{(n)}
u_{k_2,j}^{(n)} - \delta_{k_1,k_2}\right| &\leq  & \frac{C}{(\log
(1+r_n))^{1+\delta}}. \label{ii}
\end{eqnarray}
For application to periodograms, we will also need to consider pairs
$\mU^{(n)}, \mV^{(n)}$ of such matrices , and then we will assume
that in addition we have
\begin{equation}\label{iii}
\max_{1\leq k_1,k_2\leq r_n} \left|\sum_{j=1}^n u_{k_1,j}^{(n)}
v_{k_2,j}^{(n)}\right| \leq   \frac{C}{(\log (1+r_n))^{1+\delta}}.
\end{equation}
An example of a sequence of such pairs $\mU^{(n)}, \mV^{(n)}$ of
matrices with $r_n\leq \lfloor\frac{n-1}{2}\rfloor$  is
\begin{equation}\label{U+V}
u_{k,j}^{(n)}={\sqrt{\frac{2}n}}\cos\left(\frac{2\pi j k
}{n}\right),\; v_{k,j}^{(n)}=\sqrt{\frac{2}n}\sin\left(\frac{2\pi j
k }{n}\right),\; 1\leq j\leq n, 1\leq k\leq r_n.
\end{equation}
Then \eqref{i} holds trivially, while \eqref{ii}, and \eqref{iii}
follow from the fact that  $2r_n<n$ and the following well known
trigonometric identities hold for  $1\leq k_1<k_2\leq n$:
\begin{eqnarray}
%\sum_{j=1}^n\cos^2\left(\frac{2\pi jk }{n}\right)&=&
%\left\{\begin{array}{ll} n/2 &\mbox{ if $2k\ne n$}\\
%n & \mbox{ if $2k=n$}
%\end{array}
%\right.,\\
%\sum_{j=1}^n\sin^2\left(\frac{2\pi jk }{n}\right)&=&
%\left\{\begin{array}{ll} n/2 &\mbox{ if $2k\ne n$}\\
%0 & \mbox{ if $2k=n$}
%\end{array}
%\right.,\\
\label{trig 1}  \sum_{j=1}^n\cos \left(\frac{2\pi k_1 j
}{n}\right)\cos
\left(\frac{2\pi k_2 j }{n}\right)&=& \left\{\begin{array}{ll} 0 &\mbox{ if $k_1+k_2\ne n$}\\
n/2 & \mbox{ if $k_1+k_2=n$}
\end{array}
\right.,\\
 \label{trig 2} \sum_{j=1}^n\sin \left(\frac{2\pi k_1 j }{n}\right)\sin
\left(\frac{2\pi k_2 j }{n}\right)&=& \left\{\begin{array}{ll} 0 &\mbox{ if $k_1+k_2\ne n$}\\
-n/2 & \mbox{ if $k_1+k_2=n$}
\end{array}
\right.,\\
\label{trig 3} \sum_{j=1}^n\cos \left(\frac{2\pi k_1 j
}{n}\right)\sin \left(\frac{2\pi k_2 j }{n}\right)&=& 0.
\end{eqnarray}
In addition, for $1\leq k\leq n$, \begin{equation} \label{trug 4}
\sum_{j=1}^n\cos^2\left(\frac{2\pi jk
}{n}\right)=n-\sum_{j=1}^n\sin^2\left(\frac{2\pi jk }{n}\right)=
\left\{\begin{array}{ll} n/2 &\mbox{ if $2k\ne n$}\\
0 & \mbox{ if $2k=n$}
\end{array}
\right..\end{equation}

Denote by $\Phi(x)$  the standard normal cumulative distribution
function; $N(m,\sigma^2)$   denotes the normal law with mean $m$ and
variance $\sigma^2$; $X_n\Rightarrow X$ denotes weak convergence of
laws. All random variables  are assumed to be defined on a common
probability space $(\Omega,\mathcal{F}, P)$.

To avoid cumbersome notation, we state the theorem for the
univariate and bivariate cases only.  The $m$-variate extension
requires introducing $m$ sequences of matrices that satisfy
conditions \eqref{i} and \eqref{ii}, with each pair from the $n$-th
level satisfying condition \eqref{iii}; the proof requires only
minor changes.
\begin{theorem}[Almost sure CLT]\label{T1}
Suppose $X_1, X_2, \dots$  are real independent random variables
such that $E(X_j)=0$, $E(X_j^2)=1$, and $\sup_kE|X_k|^3<\infty$. Fix
$r_n\nearrow \infty$.

\begin{enumerate}
\item For $n=1,2,\dots$, let $\mU^{(n)}$ be an $r_n\times n$ matrix
such that conditions \eqref{i}, \eqref{ii} hold with some $C,\delta$
independent of $n$. Define
\begin{equation}\label{S}
S_{n,k}=\sum_{j=1}^n u_{k,j}^{(n)} X_j,\;
k=1,2,\dots,r_n.\end{equation} Then with probability $1$,
\begin{equation}\label{Conclusion1}
 \lim_{n\to\infty}  \sup_{x\in\RR} \left|r_n^{-1}\sum_{k=1}^{r_n} 1_{\{S_{n,k}\leq x\}} -\Phi(x)\right|=0.
\end{equation}
\item For $n=1,2,\dots$, let $\mV^{(n)}$ be
an $r_n\times n$ matrix such that conditions \eqref{i}, \eqref{ii}
hold for $v_{k,j}^{(n)}$ in place of $u_{k,j}^{(n)}$, and such that
the pairs $\left( \mU^{(n)}, \mV^{(n)}\right)_{n=1,2\dots,}$ satisfy
condition \eqref{iii}. Define
\begin{equation}\label{T}
 T_{n,k}=\sum_{j=1}^n
v_{k,j}^{(n)} X_j,\; k=1,2,\dots,r_n.\end{equation} Then there is a
measurable set $\Omega_0\subset\Omega$ of probability 1 such that
for all $x,y\in\RR$,
\begin{equation}\label{Conclusion}
 \lim_{n\to\infty}  r_n^{-1}\sum_{k=1}^{r_n} 1_{\{S_{n,k}\leq x,
T_{n,k}\leq y\}} =\Phi(x)\Phi(y) \mbox{ on $\Omega_0$}.
\end{equation}
\end{enumerate}
\end{theorem}
We also have the companion weak limit theorem and the large
deviation principle for the univariate case with trigonometric
coefficients given by \eqref{U+V} under the restrictions on the rate
of growth of $r_n$. The most interesting case,  $r_n=[(n-1)/2]$,
which corresponds to the spectral measures of  random circulant
matrices, is unfortunately not covered by our result;  LDP for
spectra of other random matrices are in \cite[Chapter
5]{Hiai-Petz00}.

\begin{proposition}\label{T2}
Suppose $X_1, X_2, \dots$  are real independent random variables
such that $E(X_j)=0$, $E(X_j^2)=1$, and for some constant $\tau>0$,
\begin{equation}\label{Sakhanenko}
\sup_kE\left(|X_k|^3\exp(|X_k|/\tau)\right)\leq \tau. \end{equation}
 Let
$\mU^{(n)}$ be given by \eqref{U+V}, and suppose that $r_n\to\infty$
is such that
\begin{equation}\label{r-UB}
\frac{(\log n)^2}{n} r_n^3\to 0.
\end{equation}
Then for  $x\in\RR$,
\begin{equation}
  \label{CLT}
  \frac{1}{\sqrt{r_n}}\sum_{k=1}^{r_n}\left(1_{S_{n,k}\leq x} -
  \Phi(x)\right)\Rightarrow N\left(0,\Phi(x)(1-\Phi(x))\right) \mbox{ as $n\to\infty$,}
\end{equation}

\end{proposition}
Condition \eqref{Sakhanenko}  holds for the i.i.d. sequences with
$\tau$ that depends on the law of $X_1$ when $E(\exp(\delta
|X|))<\infty$ for some $\delta>0$. We remark that \eqref{CLT} holds
true with assumptions \eqref{Sakhanenko} and \eqref{r-UB} replaced
by the assumption that there is $p>0$ such that
$$\sup_kE(|X_k|^{2+p})<\infty, \; \mbox{ and } r_n^3n^{-p/(2+p)} \to 0.$$
This follows from our proof, after substituting \cite[Section 5,
Corollary 5]{Sakhanenko-91} for  Lemma \ref{LS4}.

The large deviation principle (LDP) was motivated by \cite[Theorem
1]{March-Seppalainen-97}, which gives the LDP  from the
Brosamler-Schatte almost sure CLT. To formulate the result we need
more notation. Let $\mathcal{M}_1(\RR)$ denote the Polish space of
probability measures on the Borel sets of $\RR$ with the topology of
weak convergence. For $S_{n,k}$ given by \eqref{S}, consider the
empirical measures
\begin{equation}\label{++}
  \mu_n=\frac{1}{r_n}\sum_{k=1}^{r_n}\delta_{S_{n,k}}
\end{equation}
The rate function $I:\mathcal{M}_1(\RR)\to[0,\infty]$ in our LDP is
the relative entropy of with respect to the standard normal law,
i.e. if $\phi(x)$ denotes the normal density and
$\nu\in\mathcal{M}_1(\RR)$, then $I(\nu)=\int \ln
\frac{f(x)}{\phi(x)} f(x) dx$ if $\nu$ has the density $f$ with
respect to the Lebesgue measure and the expression is integrable,
and $I(\nu)=\infty$ otherwise. It is well known that the level sets
$I^{-1}[0,a]$ are compact for $a<\infty$.

The conclusion of the next result is the LDP  of the empirical
measures $\mu_n$ with speed $r_n$ and the rate function $I(\cdot)$.
\begin{proposition}\label{T3}
Suppose matrices $\mU^{(n)}$ and random variables $X_1,X_2,\dots$
are as in Proposition \ref{T2}. If $r_n$ is such that
$$\frac{r_n^4}{n}\to 0 \mbox{ and } \frac{\log n}{r_n}\to 0,$$
then  for all open sets $G$ and closed sets $F$ in $\mathcal
{M}_1(\RR)$,
$$
\liminf_{n\to\infty} \frac{1}{r_n}\ln\Pr(\mu_n\in G)\geq
-\inf_{\nu\in G}I(\nu)
$$
and
$$
\limsup_{n\to\infty} \frac{1}{r_n}\ln\Pr(\mu_n\in F)\leq
-\inf_{\nu\in F}I(\nu).
$$
\end{proposition}

 We   postpone the proofs to section \ref{Section: Proof}, and
we first give  some applications.

\subsection{Application to periodograms}
The periodogram of a sequence $(X_j)$ is
$$
I_n(\la)=\frac{1}{\sqrt{n}}\left|\sum_{j=1}^n  e^{-ij\la
}X_j\right|^2.
$$
The empirical distribution of the periodogram is the (random) CDF
$$F_n(x):=\frac{1}{r_n}\sum_{k=1}^{r_n} 1_{\{I_n(2\pi k/n)\leq x\}}, \; r_n=\lfloor \frac{n-1}{2}\rfloor.$$
Theorem \ref{T1} strengthens the conclusion of \cite[Proposition
4.1]{Kokoszka-Mikosch-00}   to almost sure convergence at the
expense of the assumption that third moments are finite.
\begin{corollary}\label{C1} If $X_1,X_2,\dots$ are independent with $E(X_j)=0$,
$E(X_j^2)=1$, $E(|X_j|^3)\leq M<\infty$, then $\sup_{x\geq 0}|
F_n(x)-(1-e^{-x})|\to 0$ with probability 1.
\end{corollary}
\begin{proof} Following \cite[(2.1)]{Kokoszka-Mikosch-00}, we write
$F_n(x)=\mu_n((-\infty,x])$ as the CDF of  the empirical measure
$$
\mu_n=\frac{1}{r_n}\sum_{k=1}^{r_n} \delta_{S_{n,k}^2+T_{n,k}^2},
$$
where $S_{n,k}$ and $T_{n,k}$ are defined by \eqref{S} and \eqref{T}
with
 matrices $\mU^{(n)},\mV^{(n)}$ given by \eqref{U+V}. The
result follows from Theorem \ref{T1}: if $h:E\to F$ is a continuous
mapping of Polish spaces and discrete measures
$\frac{1}{n}\sum_{k=1}^n \delta_{x_k}$ converge weakly to some
probability measure $\nu$ on the Borel sigma-field of  $E$, then the
discrete measures $\frac{1}{n}\sum_{k=1}^n \delta _{h(x_k)}$
converge weakly to the probability measure $\nu\circ h^{-1}$, see
e.g. \cite[Theorem 29.2]{Billingsley-95}. We apply this to
$h:\RR^2\to\RR$ given by $h(x,y)=x^2+y^2$ and to $\nu$  on Borel
sets of $\RR^2$, which is the product of the standard normal laws.
Then $\nu\circ h^{-1}$ is the Chi-Squared law with $2$ degrees of
freedom, i.e. it is the standard exponential law with the CDF given
by $1-e^{-x}$ for $x\geq 0$. Since the limit
$\left(1-e^{-x}\right)_+$ is a continuous CDF, it is well known, see
 \cite[Exercise 14.8]{Billingsley-95}, that the convergence is
uniform with respect to $x$.
\end{proof}
\subsection{Application to symmetric circulant and reverse circulant matrices}
 \begin{corollary} The weak
convergence in \cite {Bose-Mitra-02} and \cite[Theorem
5]{Bose-Chatterjee-Gangopadhyay-03} holds with probability one.
\end{corollary}
\begin{proof} Ref.
\cite{Bose-Mitra-02} and \cite[Theorem
5]{Bose-Chatterjee-Gangopadhyay-03} analyze the asymptotic spectrum
of the $n\times n$ symmetric random matrices with the typical
eigenvalues of the form $\pm \sqrt{S_{n,k}^2+T_{n,k}^2}$, where
$k=1,2,\dots, \lfloor\frac{n-1}{2}\rfloor $, and
$\mU^{(n)},\mV^{(n)}$ defined by \eqref{U+V}, see \cite[Lemma
1]{Bose-Mitra-02}. Omitting at most two eigenvalues does not change
the convergence of the spectral measure, so theorem \ref{T1} implies
that the convergence holds with probability one by the argument
similar to the proof of Corollary \ref{C1}.
\end{proof}
Suppose $\mA_n$ is a symmetric random circulant matrix formed from
the independent random  variables by taking as  the first row
$[\mA_n]_{1,j} =X_j$, $j=1,2,\dots,[(n+1)/2]$ and   $[\mA_n]_{1,j}
=[\mA_n]_{1,n-j}$ for other $j$. The next corollary strengthens
\cite[Remark 2]{Bose-Mitra-02} to almost sure convergence, and
relaxes the integrability and i.i.d assumption in \cite[Theorem
1.5]{Miller05}. To justify the later claim, we note
 that a "palindromic matrix" analyzed in \cite{Miller05} differs from
 $\mA_n$ by the last row and column only; thus their ranks differ  by
 at most  one, and asymptotically "palindromic matrices" and random circulant matrices have
 the same spectrum, see \cite[Lemma 2.2]{Bai99}.
 \begin{corollary} If $X_1, X_2,\dots$ are independent with uniformly bounded third moments, with common mean $E(X_j)=m$
 and common variance $Var(X_j)=\sigma^2>0$, then the
 spectrum of
 $\frac{1}{\sigma\sqrt{n}}\mA_n$ converges weakly with probability 1 to
 the standard normal law.
 \end{corollary}
 \begin{proof} Subtracting the rank 1 matrix does not change the asymptotic of the spectrum, thus without loss of generality
 we may assume $m=0$; rescaling the variables by $\sigma>0$ we can assume $E(X_j^2)=1$.
With the exception of at most two eigenvalues, the remaining
eigenvalues of $\mA_n/\sqrt{n}$ are of multiplicity two and are
given by \eqref{S}
 with the trigonometric matrix
 $\mU^{(n)}$ given by \eqref{U+V}, see \cite[Remark
 2]{Bose-Mitra-02}. Thus the weak convergence with probability one of the spectral law of $\mA_n/\sqrt{n}$
 to $N(0,1)$ follows from Theorem
 \ref{T1}.
 \end{proof}

\subsection{Application to random orthogonal matrices}
A well known result of Poincar\'e says that if $\mU^{(n)}$ is a
random orthogonal matrix uniformly distributed on $O(n)$ and
$x_n\in\RR^n$ is a sequence of vectors of norm $\sqrt{n}$ then the
first $k$ coordinates of  $\mU^{(n)} x_n$ are asymptotically normal
and independent, see e.g. \cite[Exercise 29.9]{Billingsley-95}.
\begin{corollary}\label{P-Haar} Let  $X_1,X_2,\dots$ be independent, $E(X_j)=0$,
$E(X_j^2)=1$ and $\sup_jE(|X_j|^3)<\infty$. Let $\mU^{(n)}$ be a
random orthogonal matrix uniformly distributed on $O(n)$ and
independent of $(X_j)$. Define $S_{n,k}$ by \eqref{S}. Then with
probability one
\begin{equation}
  \label{Haar-CLT}  \sup_{x\in\RR}\left|\frac{1}{n}\sum_{k=1}^n1_{\{S_{n,k}\leq x\}} -
\frac{1}{\sqrt{2\pi}}\int_{-\infty}^x e^{-u^2/2}du\right|\to 0 \;
\mbox{ as $n\to \infty$.}
\end{equation}
\end{corollary}
\begin{proof} This result has a direct elementary proof, which we
learned from Jack Silverstein. This proof shows that the result
holds true also for i.i.d. random variables with finite second
moments. Here we derive it as a corollary to Theorem \ref{T1}.

Orthogonal matrices satisfy \eqref{ii} with $r_n=n$. By
\cite[Theorem 1]{Jiang-05},  \eqref{i} holds with probability $1$.
Therefore, redefining  $\mU^{(n)}$ and $(X_j)$ on the product
probability space $\Omega_U\times\Omega_X$, by \cite[Theorem
1]{Jiang-05}, there is a subset $\Omega_U'$ of probability $1$ such
that for each $\omega_1\in\Omega_U'$ by   Theorem \ref{T1} there is
a measurable subset $\Omega_{X,\omega_1}\subset \Omega_X$ of
probability one such that \eqref{Haar-CLT} holds. By Fubini's
Theorem, the set of all pairs $(\omega_1,\omega_2)$ for which
\eqref{Haar-CLT} holds has probability one.
\end{proof}
\section{Proofs}\label{Section: Proof}
\subsection{Proof of Theorem \protect{\ref{T1}}}
The proof of   consists of several lemmas.

\begin{lemma}[{\cite[Theorem 2.6]{berti-pratelli-rigo-06}}]\label{L1} It is enough to verify almost sure convergence of
characteristic functions.
\end{lemma}

To use Lemma \ref{L1}, fix   real $s,t$ and consider the (random)
characteristic function

$$\Phi_n(s,t,\omega)=r_n^{-1}\sum_{k=1}^{r_n} \exp(is S_{n,k}(\omega)+i t T_{n,k}(\omega)).$$

%Lemma 2. $E(\Phi_n(t,\omega) -> e^{-t^2/2}
%
%This  essentially follows from the standard proof of the CLT. The
%error per each term under the sum  is estimated by C |t|^3 \max_j
%|u_{k,j}^{(n)}(n)|, so the convergence follows from the uniform
%smallness assumption (i).
%
%It seems I do not need this lemma anyway. Instead I directly proceed
%with.

\begin{lemma} \label{L3}   There is $C=C(t)$ that does not depend on $n>1$ such
that for all large enough $n>N(t)$ we have

\begin{equation}\label{*}
E| \Phi_n(s,t,\omega) -  e^{-s^2/2-t^2/2}|^2 \leq C/(\log (1+r_n)
)^{1+\delta}.
\end{equation}
\end{lemma}
\begin{proof}
The left hand side of \eqref{*} is
$$
\frac{1}{r_n^2}\sum_{k_1,k_2=1}^{r_n}E\left(\left(e^{i s S_{n,k_1}+i
t T_{n,k_1}}-e^{-(s^2+t^2)/2}\right)\left(e^{-i s S_{n,k_2}-i t
T_{n,k_2}}-e^{-(s^2+t^2)/2}\right)\right).
$$
Denote $\varphi_j(t)=E(e^{i t X_j})$. Then  the left hand side of
\eqref{*} can be bounded by
\begin{multline*}
  2/r_n+\frac{1}{r_n^2}\sum_{k_1\ne k_2} \Big|\prod_{j=1}^n\varphi_j\left(s(u_{k_1,j}^{(n)}-u_{k_2,j}^{(n)})+t(v_{k_1,j}^{(n)}-v_{k_2,j}^{(n)})\right)+e^{-t^2-s^2}\\
  -e^{-(s^2+t^2)/2} \prod_{j=1}^n\varphi_j(su_{k_1,j}^{(n)}+tv_{k_1,j}^{(n)})- e^{-(s^2+t^2)/2}
  \prod_{j=1}^n\varphi_j(-su_{k_2,j}^{(n)}-tv_{k_2,j}^{(n)})\Big|\\
\leq
2/r_n+\max_{k_1\ne k_2}\left|\prod_{j=1}^n\varphi_j\left(s(u_{k_1,j}^{(n)}-u_{k_2,j}^{(n)})+t(v_{k_1,j}^{(n)}-v_{k_2,j}^{(n)})\right)-e^{-t^2-s^2}\right|\\
+2\max_{ \pm}\max_{1\leq k\leq r_n}\left| \prod_{j=1}^n\varphi_j(\pm
su_{k ,j}^{(n)}\pm tv_{k ,j}^{(n)})- e^{-(s^2+t^2)/2}\right|.
\end{multline*}
We will show how to bound the middle term, as the last one is
handled similarly. Trivially,
 \begin{multline} \max_{k_1\ne
k_2}\left|\prod_{j=1}^n\varphi_j\left(s(u_{k_1,j}^{(n)}-u_{k_2,j}^{(n)})+t(v_{k_1,j}^{(n)}-v_{k_2,j}^{(n)})\right)-e^{-t^2-s^2}\right|
\leq A_n+B_n,
\end{multline}
where
\begin{multline}
A_n=\max_{k_1\ne
k_2}\Big|\prod_{j=1}^n\varphi_j\left(s(u_{k_1,j}^{(n)}-u_{k_2,j}^{(n)})+t(v_{k_1,j}^{(n)}-v_{k_2,j}^{(n)})\right)\\
-
\prod_{j=1}^n\exp\left(-\frac12\left(s(u_{k_1,j}^{(n)}-u_{k_2,j}^{(n)})+t(v_{k_1,j}^{(n)}-v_{k_2,j}^{(n)})\right)^2\right)\Big|,
\end{multline}
\begin{multline}B_n= \max_{k_1\ne k_2}
\left|\exp\left(-\frac12\sum_{j=1}^n\left(s(u_{k_1,j}^{(n)}-u_{k_2,j}^{(n)})+t(v_{k_1,j}^{(n)}-v_{k_2,j}^{(n)})\right)^2\right)-e^{-t^2-s^2}\right|.
\end{multline}
By \eqref{i}, for large enough $n$   we have \begin{equation}
\label{**} 0\leq
\max_{j,k_1,k_2}\left|s(u_{k_1,j}^{(n)}-u_{k_2,j}^{(n)})+t(v_{k_1,j}^{(n)}-v_{k_2,j}^{(n)})\right|\leq
1.
\end{equation} We now use the well known bound, see e.g.
\cite[(27.13)]{Billingsley-95},
 which bounds the first term by
 \begin{multline*}
A_n\leq  C\max_{k_1, k_2}\sum_{j=1}^n
 \left|s(u_{k_1,j}^{(n)}-u_{k_2,j}^{(n)})+t(v_{k_1,j}^{(n)}-v_{k_2,j}^{(n)})\right|^3 \\ \leq
 \frac{C(s,t)}{(\log (1+r_n))^{1+\delta}} \left(\max_{k}\sum_{j=1}^n
 (u_{k,j}^{(n)})^2+\max_{k}\sum_{j=1}^n(v_{k,j}^{(n)})^2\right)\leq \frac{C'(s,t)}{(\log
(1+r_n))^{1+\delta}}.
 \end{multline*}
The second term is bounded using \eqref{ii} and \eqref{iii} as
follows. We first note that   $k_1\ne k_2$ implies
$$
\sum_{j=1}^n\left(s(u_{k_1,j}^{(n)}-u_{k_2,j}^{(n)})+t(v_{k_1,j}^{(n)}-v_{k_2,j}^{(n)})\right)^2=2s^2+2t^2+
R(n,k_1,k_2,s,t),$$ where
\begin{multline*}
|R(n,k_1,k_2,s,t)| \leq
s^2\left|\sum_{j=1}^n\left((u_{k_1,j}^{(n)})^2+(u_{k_2,j}^{(n)})^2\right)-2\right|
+t^2\left|\sum_{j=1}^n\left((v_{k_1,j}^{(n)})^2+(v_{k_2,j}^{(n)})^2\right)-2\right|
\\ +2s^2 \left|\sum_{j=1}^n u_{k_1,j}^{(n)}u_{k_2,j}^{(n)}\right| +
2t^2 \left|\sum_{j=1}^n v_{k_1,j}^{(n)}v_{k_2,j}^{(n)}\right|\\+2
\left|st\sum_{j=1}^n
u_{k_1,j}^{(n)}v_{k_1,j}^{(n)}-u_{k_2,j}^{(n)}v_{k_1,j}^{(n)}+u_{k_2,j}^{(n)}v_{k_2,j}^{(n)}-u_{k_1,j}^{(n)}v_{k_2,j}^{(n)}\right|
 \leq \frac{C(s,t)}{(\log (1+r_n))^{1+\delta}}.
\end{multline*}
Since $|e^{a}-e^{b}|\leq   |a-b|\max_{u\in[a,b]}e^u$, and \eqref{**}
holds, therefore for large enough $n$ this implies
%\begin{multline*}
$$
%\max_{k_1\ne k_2}
%\left|\exp\left(-\frac12\sum_{j=1}^n\left(s(u_{k_1,j}^{(n)}-u_{k_2,j}^{(n)})+t(v_{k_1,j}-v_{k_2,j}^{(n)})\right)^2\right)-e^{-t^2-s^2}\right|
B_n\leq \frac{C(s,t)}{(\log (1+r_n))^{1+\delta}}. $$
%\end{multline*}
\end{proof}

To prove almost sure convergence we will use the following.

\begin{lemma}[{\cite[Theorem 1]{Lyons-88}}]\label{L4} Let   $Y_1, Y_2, \dots$ be uniformly bounded $\CC$-valued and possibly dependent
random variables, $r_n\nearrow\infty$. Suppose
$Z_n=\frac1{r_n}\sum_{k=1}^{r_n} Y_k$. If $E|Z_n|^2\leq C/ (\log
(1+r_n))^{1+\delta}$ then $Z_n \to 0$ with probability 1.
\end{lemma}
%\begin{proof}
%For completeness, we give a short proof suggested by \cite[Exercise
%6.7]{Billingsley-95}. Given $q>1$, let $N_k=\min\{n: r_n\geq q^k\}$.
%Then $Z_{N_k}\to 0$ with probability $1$ by the Borel-Cantelli
%lemma.
%
%Now given $n\geq 1$, choose $k=k(n)$ such that $N_{k-1}\leq n <
%N_k$. Then \begin{multline*}  |Z_n|\leq |Z_{N_{k-1}}|+
%\frac1{r_{N_{k-1}}}\max_{ N_{k-1}\leq n < N_k}
%\sum_{j=r_{N_{k-1}}}^{r_n} |Y_j| \leq |Z_{N_{k-1}}|+ C (q-1).
%\end{multline*}  Since $q>1$ can be chosen arbitrarily close to $1$,
%this shows that $\limsup_n |Z_n| =0$ on the same subset of $\Omega$
%on which $Z_{N_k}\to 0$.
%\end{proof}
\begin{proof}[Proof of Theorem \protect{\ref{T1}}]

To prove part (ii), let $Y_k =\exp(i s S_{n,k}+i t
T_{n,k})-e^{-s^2/2-t^2/2}$. Then by Lemma \ref{L3}, the assumptions
of  Lemma \ref{L4} are satisfied. So
$Z_n(\omega)=\Phi_n(s,t,\omega)-e^{-s^2/2-t^2/2}\to 0$ with
probability one. By Lemma \ref{L1}, this implies \eqref{Conclusion}.

The proof part (i) is similar, and essentially consists of taking
$t=0$ in the above calculations; once we establish the weak
convergence on a set $\Omega_0$ of probability $1$, due to
continuity of $\Phi(x)$, the convergence is uniform in $x$ for every
$\omega\in\Omega_0$, see \cite[Exercise 14.8]{Billingsley-95}.

\end{proof}

\subsection{Proof of Propositions \ref{T2} and \ref{T3}}

The proofs rely on strong approximation of the partial sum processes
indexed by the Lipschitz functions $f_k(x)=\cos(2\pi k x)$, compare
\cite[Theorems 2.1, 2.2]{Grama-Nussbaum-02}. We derive suitable
approximation directly from the following result.

\begin{lemma}[Sakhanenko {\cite[Theorem 1]{Sakhanenko-84}}]\label{LS4} Suppose
$X_1,X_2,\dots$ satisfy the assumptions of Proposition \ref{T2}.
Then there is a constant $c$ such that for every $n$ one can realize
$X_1,X_2,\dots,X_n$ on a probability on which there are i.i.d.
$N(0,1)$ random variables
$\widetilde{X}_1,\widetilde{X}_2,\dots,\widetilde{X}_n$ such that
the partial sums $S_j=\sum_{i=1}^jX_i$ and
$\widetilde{S}_j=\sum_{i=1}^j\widetilde{X}_i$ satisfy
\begin{equation}
  \label{Exp bound}
  E \left(\exp\left(\frac{c}{\tau} \max_{1\leq j \leq n}
  |S_j-\widetilde{S}_j|\right)\right)\leq 1+\frac{n}{\tau}.
\end{equation}
(Recall that $\tau$ is defined in \eqref{Sakhanenko}.)
\end{lemma}
We use this lemma as follows. For every $n$, we redefine
$X_1,X_2,\dots,X_n$ onto a new probability space
$(\Omega_n,\mathcal{F}_n,P_n)$ on which we have the i.i.d. standard
normal r.v. $\widetilde{X}_1,\widetilde{X}_2,\dots,\widetilde{X}_n$
which satisfy \eqref{Exp bound}.  We then define $\{S_{n,k}:
k=1,2,\dots,r_n\}$ by \eqref{S}, and we also define
\begin{equation}
  \label{S-tilde}
  \widetilde{S}_{n,k}=\sum_{j=1}^n \widetilde{X}_j\cos(\frac{2\pi j
  k}{n}), \; k=1,2,\dots,r_n.
\end{equation}
(Clearly we should have used the triangular array  notation
$(X_{j,n})_j$ instead of $(X_j)_j$; we can safely omit the subscript
$n$ here, since its re-appearance in the partial sums $S_{n,k}$
keeps track of $n$ anyway.)

The assumptions and the conclusions of Propositions \ref{T2} and
\ref{T3} are not affected by such a change.

We note that from the trigonometric identities (\ref{trig
1}-\ref{trug 4}), it follows that for fixed $n$ random variables
$\widetilde{S}_{n,1},\widetilde{S}_{n,2},\dots,\widetilde{S}_{n,r_n}$
are i.i.d standard normal. Therefore, if $r_n\to\infty$ then for any
$\eta$ we have
\begin{equation}
\label{T2-normal}
\frac{1}{\sqrt{r_n}}\sum_{k=1}^{r_n}\left(1_{\widetilde{S}_{n,k}\leq
x+\eta /\sqrt{r_n}} -
  \Phi\left(x+\frac{\eta}{\sqrt{ r_n}}\right)\right)\Rightarrow N\left(0,\Phi(x)(1-\Phi(x))\right) \mbox{ as $n\to\infty$.}
\end{equation}
(This is just the normal approximation to the binomial random
variables with the probabilities of success $
\Phi\left(x+\frac{\eta}{\sqrt{ r_n}}\right)\to \Phi(x)$.)

Consider now the corresponding empirical measures
\begin{equation}
  \label{T3 normal}
  \widetilde{\mu}_n=\frac{1}{r_n}\sum_{k=1}^{r_n}\delta_{\widetilde{S}_{n,k}}
\end{equation}
 If $r_n\to\infty$, then by Sanov's Theorem, see e.g. \cite{D-Z-98},
for all open sets $G$ and closed sets $F$ in $\mathcal {P}(\RR)$,
$$
\liminf_{n\to\infty} \frac{1}{r_n}\ln\Pr(\widetilde{\mu}_n\in G)\geq
-\inf_{\nu\in G}I(\nu)
$$
and
$$
\limsup_{n\to\infty} \frac{1}{r_n}\ln\Pr(\widetilde{\mu}_n\in F)\leq
-\inf_{\nu\in F}I(\nu).
$$
The plan of proof is to deduce Propositions \ref{T2} and \ref{T3}
from these two facts.

By taking a product space, without loss of generality we assume that
all random variables $S_{n,k},\widetilde{S}_{n,k}$, $n\geq 1, 1\leq
k\leq r_n$, are defined on the common probability space. We never
need joint distributions of these variables for different $n$, but
such a choice simplifies the notation.

\begin{proof}[Proof of Proposition \ref{T2}]
Denote $$ Z_n(x)=\frac{1}{\sqrt{r_n}}\sum_{k=1}^{r_n}1_{S_{n,k}\leq
x}$$ and let $\widetilde{Z}_n(x)$ denote the corresponding sum for
the i.i.d. normal random variables \eqref{S-tilde} from Lemma
\ref{LS4}. Fix $\eps>0$, and let $\eps_n=\eps\sqrt{2\pi/r_n}$. From
the trivial bound
$$\left|\Phi(x\pm\eps_n)-\Phi(x)\right|\leq \eps_n/\sqrt{2\pi}$$
we get
\begin{multline}\label{+++}
 \widetilde{Z}_n(x-\eps_n)-\sqrt{r_n}\Phi(x-\eps_n) -\eps-R_n \\
  \leq Z_n(x)-\sqrt{r_n}\Phi(x)\leq
  \widetilde{Z}_n(x+\eps_n)-\sqrt{r_n}\Phi(x+\eps_n)+\eps+R_n,
\end{multline}
where \begin{equation*} R_n=\sqrt{r_n}\,1_{A_n},\; A_n=\{\max_{1\leq
k \leq r_n} |S_{n,k}-\widetilde{S}_{n,k}|>\eps_n\}.
\end{equation*}
Since $\cos(\cdot)$ is a Lipschitz function,
 \begin{multline}\label{+}    \max_{1\leq k \leq r_n}
|S_{n,k}-\widetilde{S}_{n,k}|= \max_k
\left|\sum_{j=1}^n(X_j-\widetilde{X}_j)\cos\left(\frac{2\pi j
k}{n}\right)\right|/\sqrt{n}\\
\leq n^{-1/2}\max_{k\leq r_n}
\left|\sum_{j=1}^{n-1}(S_j-\widetilde{S}_j)\left(\cos\left(\frac{2\pi
j k}{n}\right)-\cos\left(\frac{2\pi (j+1)
k}{n}\right)\right)\right|+n^{-1/2}|S_n-\widetilde{S}_n| \\
\leq n^{-1/2}\max_{k\leq r_n}
\sum_{j=1}^{n-1}|S_j-\widetilde{S}_j|\frac{2\pi k}{n}+n^{-1/2}|S_n-\widetilde{S}_n| \\
\leq \frac{1+2\pi r_n}{\sqrt{n}}\max_{1\leq j \leq
n}|S_j-\widetilde{S}_j|.
\end{multline}
From Lemma \ref{LS4} we see that for large enough $n$ so that
$r_n/\sqrt{n}\leq c/\tau$ we have by Markov inequality
\begin{multline*}
  \Pr(A_n)\leq \Pr\left(\max_{1\leq j \leq n}|S_j-\widetilde{S}_j|\geq
  \frac{\eps\sqrt{2\pi n}}{\sqrt{r_n}(1+2\pi r_n)}
\right) \\
\leq \exp\left(-\frac{c\eps\sqrt{2\pi n}}{\tau\sqrt{r_n}(1+2\pi
r_n)} \right)E\left(\exp \left(c/\tau\max_{1\leq j \leq
n}|S_j-\widetilde{S}_j|\right)\right)\\
\leq \exp\left(\log(1+\frac{ n}{\tau}) -\frac{c\eps\sqrt{2\pi
n}}{\tau\sqrt{r_n}(1+2\pi r_n)} \right)\to 0.
\end{multline*}
Thus $R_n\to 0$ in probability. Since $\eps>0$ is arbitrary,
Proposition \ref{T2} follows from \eqref{+++} by \eqref{T2-normal}.
\end{proof}
Our proof of Proposition \ref{T3} is based on the following
approximation lemma.
\begin{lemma}[{\cite[Theorem 4.9]{Baxter-Jain-96}}] \label{L5} Suppose the families of real random
variables $S_{n,k},\widetilde{S}_{n,k}$,$n\geq 1, 1\leq k \leq r_n$
are such that for every $\theta>0$
\begin{equation}
\label{@b} \limsup _{n\to \infty } {r_n^{-1}}\log E\left(\exp
\left(\theta
\sum_{k=1}^{r_n}|S_{n,k}-\widetilde{S}_{n,k}|\right)\right) \leq 1.
\end{equation}
 If the empirical measures \eqref{T3 normal}
satisfy the LDP in $\mathcal{M}_1(\RR)$ with speed $r_n$ and a good
rate function $I(\cdot)$, then the empirical measures \eqref{++}
satisfy the LDP in $\mathcal{M}_1(\RR)$ with the same speed $r_n$
and the same rate function $I(\cdot)$.
\end{lemma}
%\begin{proof} This is a consequence of the fact that empirical measures
%$\mu_n$ and $\widetilde{\mu}_n$ are exponentially equivalent, so
%\cite[Theorem 4.2.13]{D-Z-98} applies. (For more complicated
%versions of this lemma,  see \cite[Theorem 4.9]{Baxter-Jain-96},
%\cite[Lemma 3.1]{Bryc95a}.)
%
% For completeness we enclose the
%proof.  We consider $\mathcal{M}_1(\RR)$ with the bounded Lipschitz
%metric
%\begin{equation} \label{BLip}
%\dbl (\mu,\nu)=\sup\{\int f d\mu-\int f d\nu:\|f\|_\infty +
%\|f\|_L\leq 1\},
%\end{equation}
%where $$\|f\|_\infty=\sup_x|f(x)|,  \; \|f\|_L=\sup_{x\ne y}
%\frac{|f(x)-f(y)|}{|x-y|},$$ see \cite[Section 11.3]{Dudley}.
%Suppose $\delta>0$. Fix $M>0$ and let $\theta=(M+1)/\delta$. Then
%\begin{multline*} \Pr(\dbl(\mu_n,\widehat{\mu}_n)>\delta)\leq
%\exp(-\delta \theta r_n) E\left(\exp ( \theta r_n
%\dbl(\mu_n,\widehat{\mu}_n)\right) \\
%=\exp(-\delta\theta r_n) E\left(\exp
%(\theta\sup_f\sum_{k=1}^{r_n}|f(S_{n,k})-f(\widetilde{S}_{n,k})|\right)\\
%\leq\exp(-\delta\theta r_n) E\left(\exp
%(\theta\sum_{k=1}^{r_n}|S_{n,k}-\widetilde{S}_{n,k}|\right).
%\end{multline*}
%By \eqref{@b}, $E\left(\exp
%(\theta\sum_{k=1}^{r_n}|S_{n,k}-\widetilde{S}_{n,k}|\right)\leq
%e^{r_n}$ for all large enough $n$. Therefore, $$\frac{1}{r_n}\ln
%\Pr(\dbl(\mu_n,\widehat{\mu}_n)>\delta)\leq -\delta\theta+1=-M.$$
%Since $M>0$ was arbitrary,
%$$\limsup_{n\to\infty}\frac{1}{r_n}\ln
%\Pr(\dbl(\mu_n,\widehat{\mu}_n)>\delta)=-\infty.$$ Thus empirical
%measures $\mu_n$ and $\widetilde{\mu}_n$ are exponentially
%equivalent, and LDP follows from \cite[Theorem 4.2.13]{D-Z-98}.
%\end{proof}
\begin{proof}[Proof of Proposition \ref{T3}]
Since the Large Deviation Principle for $\widetilde{\mu}_n$ follows
from Sanov's Theorem, to end the proof we only need to verify
assumption \eqref{@b} of Lemma \ref{L5}. Using \eqref{+}, we see
that for large $n$ there is $C$ such that
\begin{multline*}
  E\left(\exp (\theta
\sum_{k=1}^{r_n}|S_{n,k}-\widetilde{S}_{n,k}|)\right)\leq
E\left(\exp (\theta r_n \max_{1\leq k\leq
r_n}|S_{n,k}-\widetilde{S}_{n,k}|)\right) \\
\leq E\left(\exp (C r_n^2n^{-1/2} \max_{1\leq j \leq
n}|S_{j}-\widetilde{S}_{j}|)\right).
\end{multline*}
Since $r_n^2n^{-1/2}\to 0$, therefore for large enough $n$ we can
apply \eqref{Exp bound}. We get
$${r_n^{-1}}\log E\left(\exp (\theta
\sum_{k=1}^{r_n}|S_{n,k}-\widetilde{S}_{n,k}|)\right)\leq C
\frac{\log (1+n/\tau)}{r_n}\to 0.$$
\end{proof}

{\bf Acknowledgement} The second-named author (WB) thanks S. Miller
for the early version of  \cite{Miller05}, and M. Peligrad for a
helpful discussion. The authors thank J. Silverstein for showing us
the elementary direct proof of Corollary \ref{P-Haar}.
%BibTeX
%\bibliographystyle{alpha}
%\bibliographystyle{apalike}
%\bibliographystyle{newapa} %name (year)
%\bibliographystyle{plain} %full fist names
%\bibliographystyle{plain}
%\bibliographystyle{unsrt}
\bibliographystyle{acm} %initials after name
\bibliography{../Vita,asclt}
\end{document}